\documentclass[12pt]{article}
\usepackage{amssymb} 
\usepackage{amsmath}
\newtheorem{teo}{Theorem}

\newtheorem{rem}{Remark}
\newtheorem{lem}{Lemma}

\title{ Approximation by refinement masks
}
\author{ Elena A. Lebedeva
\footnote{Mathematics and Mechanics Faculty, Saint Petersburg State University,
Universitetsky prospekt, 28, Peterhof,  Saint Petersburg, 198504, Russia
 }
}
\date{
 ealebedeva2004@gmail.com
}

\begin{document}
\maketitle

\newcommand{\nul}{{\bf0}}
\newcommand{\rd}{{\mathbb R}^d}
\newcommand{\zd}{{\mathbb Z}^{d}}
\renewcommand{\r}{{\mathbb R}}
\newcommand{\z} {{\mathbb Z}}
\newcommand{\cn} {{\mathbb C}}
\newcommand{\n} {{\mathbb N}}

\begin{abstract} In the paper we design a Parseval wavelet frame  with a compact support.  The corresponding refinement mask uniformly approximates an arbitrary continuous periodic function $f$, $f(0)=1$, $|f(x)|^2+|f(x+\pi)|^2\le 1$. The refinable function has stable integer shifts. 
              
\end{abstract}



In the paper we study density of the set of all refinement masks for framelets. We design a Parseval wavelet frame with a compact support. The corresponding refinable function has stable integer shifts. And the main contribution to the topic is the following. The corresponding refinement mask uniformly approximates an arbitrary continuous periodic function $f$ such that $|f(x)|^2+|f(x+\pi)|^2\le 1$. The last restriction is natural, since any  refinement mask $m_0$ has  to satisfy the same inequality
\begin{equation}
|m_0(x)|^2+|m_0(x+\pi)|^2\le 1.
\label{main_m}
\end{equation}
 We  provide more details addressed this question later.
Besides the pure theoretical interest in properties of  wavelet frames, we are motivated by a closed connection between refinement masks and filter banks. It is well-known that the sequence of Fourier coefficients of a refinement mask forms a low-pass filter, and inequality (\ref{main_m}) provides a necessary condition for a filter bank to have the prefect reconstruction property (see \cite{BGG}). Actually the density of the set of all refinement masks    
means that we can approximate any low-pass filter with the prefect reconstruction property by a low-pass filter of compactly supported Parseval wavelet frame generating by multiresolution analysis (MRA). 

 In \cite{Le}, an arbitrary $1$-periodic continuous function $f$, $f(0)=1$, is approximated   by a one-parametric family of functions, depending on refinement masks. The results of this article   improves \cite{Le} in two directions. First, we approximate a function by a refinement mask itself. Second, a trigonometric polynomial served as a refinement mask is provided by interpolation linear summation method based on initial data.    

The density  of the set of all wavelet generators is studied in \cite{Larson}, \cite{Bownik}, \cite{CM}. It the last two papers it is proved  that  the set of all framelets is dense in $L_2(\mathbb{R})$ and the subset of tight framelets is not dense. It serves for us as a motivation as well.  
In the paper we ask the same question for refinement masks.

We recall briefly notations and results we need and draw the reader's attention to the corresponding properties of refinement masks. 
  To design wavelet frames we exploit the unitary extension principle by Ron and Shen \cite{RS}. We recall it here. In the sequel, we denote by $C$ and $L_{p}$ the spaces of $2\pi$-periodic continuous functions and $2\pi$-periodic  functions $f$ such that $|f|^p$ is integrable, respectively.  $\mathbb{T} = [-\pi,\,\pi)$   denotes the one-dimensional torus. $\delta_{k,n} = 0$ for $n \neq k,$ $\delta_{k,k}=1.$

\textbf{The unitary  extension principle   }(UEP): 
	Suppose there exist functions  $\varphi \in L_2(\mathbb{R}),$ (the refinable function) and  $m_0 \in L_{2}$  (the refinement mask) such that 
	\begin{equation}
	\label{refeq}
	\widehat{\varphi}(2 \xi) = m_0(\xi) \widehat{\varphi}(\xi), 
	\end{equation}
	and $\lim_{\xi \to 0}\widehat{\varphi}(\xi) =1.$ 
	We fix $m_1,\dots, m_q \in L_{2}$ (the wavelet masks) and define  $\psi_1,\dots,\psi_q\in L_2(\mathbb{R})$ (the wavelet generators) by 
	\begin{equation}
		\label{psi}
	\widehat{\psi_r}(2 \xi) = m_r(\xi) \widehat{\varphi}(\xi), \ \ r=1,\dots,q.
\end{equation}
	If for a.e. $\xi \in \mathbb{T}$
	\begin{equation}
		\label{matr}
		\left\{
	\begin{array}{l}
		\sum\limits_{r=0}^{q} |m_r(\xi)|^2 = 1, \\
		\sum\limits_{r=0}^{q} m_r(\xi)\overline{m_r(\xi+\pi)} = 0,
	\end{array}
	\right.
\end{equation}
	then the functions $\displaystyle \left\{\psi_{r,j,k}\right\}_{j,k\in \mathbb{Z},r=1,\dots,q}$ form a Parseval frame for $L_2(\mathbb{R}).$ 
	Here we exploit the usual notation $\psi_{r,j,k}(x):=2^{j/2}\psi_r(2^j x+k).$
	
	Revisiting our motivation, it is necessary to note that the inequality (\ref{main_m}) is a straight corollary of the system (\ref{matr}). So, it is not possible  to approximate by refinement masks a $2\pi$-periodic function $f$ which does not satisfy the inequality $\left|f(\xi)\right|^2+\left| f(\xi+\pi)\right|^2 \le 1.$  At the same time, the system (\ref{matr}) provides a necessary and sufficient condition for a filter bank to have the prefect reconstruction property (see \cite{BGG}). It adjusts ``filter'' part of our motivation to study density of masks.

There exist  refinable functions $\varphi$ such that they do not form an MRA.
It is known (see, \cite[Theorem 1.2.14]{NPS}) that the refinement equation (\ref{refeq}), continuity of $\widehat{\varphi}$ at the origin, and the stability of integer shifts of $\varphi$ (see \cite[Definition 3.4.11]{NPS}) are sufficient for a refinable function $\varphi$ to generate an MRA. The stability of shifts can be provided by specific properties of refinement mask.   
 We recall here necessary notions.  Let $g(\xi)$ be a $2\pi$-periodic function. If for $\alpha\in \mathbb{T}$ we get 
$g(\alpha)= g(\alpha+\pi)=0,$ then the pair $\{\alpha,\,\alpha+\pi\}$ is called a pair of symmetric roots of $g(\xi).$ A set of different complex numbers $\{b_1,\dots,b_n\}$ is called  cyclic if $b_{j+1}=b_j^2,$ $j=1,\dots,n$, and $b_n^2=b_1.$ In other words, $\displaystyle b_j={\rm e}^{i \beta_j}$, where  $\beta_j=2\pi m 2^{j-1}(2^n-1)^{-1}$ for some $m\in\mathbb{Z}.$ A cyclic set is called a cycle of a function $g$ if  $g(\beta_j+\pi)=0$ for all $j=1,\dots,n.$ The cycle $\{1\}$ is called trivial, all other cycles are called nontrivial.
 To provide the stability of the integer shifts we use the following statement \cite[Corollary 3.4.15]{NPS}: Integer shifts of a refinable function are stable if and only if the mask $m_0$ has neither nontrivial cycles nor a pair of  symmetric roots on $\mathbb{T}$.

\begin{rem}
\label{rem1}
 Let $m_0$ be a trigonometric polynomial, $m_0(0)=1$, and the inequality (\ref{main_m}) is fulfilled. 
 These assumptions are sufficient for  the function
$m_0$ to be  a refinement mask of the UEP. Indeed, we assume $\displaystyle \widehat{\varphi}(\xi) = \prod_{j=1}^{\infty} m_0(\xi/2^j)$.  By the Mallat theorem \cite[Lemma 4.1.3]{NPS})  $\varphi \in L_2(\mathbb{R})$, and since $\widehat{\varphi}$ is an entire function of exponential  type, it follows that  $\widehat{\varphi}$ is continuous at zero and  
$\displaystyle \widehat{\varphi}(0) = m_0(0)=1$.
Therefore, $\varphi$ is the refinement function corresponding to the refinement mask $m_0(\xi).$
 \end{rem}

Now we are ready to prove the main result of the paper. 
 
\begin{teo}
Suppose $f\in C$, $f(0)=1,$ $|f(x)|^2+|f(x+\pi)|^2\le 1$,  $\varepsilon >0$. Then there exists a compactly supported Parseval wavelet frame with a refinement mask $m_0$ such that $\|f - m_0\|_{C}<\varepsilon.$  The refinable function $\varphi$ has stable integer shifts.
\label{main}
\end{teo} 
 
\textbf{Proof}. 1.  
We approximate the function $f$ by a piecewise linear function with only finite number of roots, without nontrivial cycles and without pairs of symmetric roots. To this end, we find a piecewise linear function $f_1$ such that $\|f-f_1\|<\varepsilon/6.$ 
For example, $f_1$ interpolates $f$ at the equidistant points 
$\{\xi_k\}_{k=0}^{2n-1} =\{\pi k/n\}_{k=0}^{2n-1}.$  
Now we change the values of $f_1$ in the neighborhoods of the  segments, where $f_1\equiv 0$, if any. 
Let $f(\xi_k)=0$ for $k=i,\dots,j,$ and $f(\xi_{i-1}) \neq 0,$ $f(\xi_{j+1}) \neq 0.$ There are three possibilities. 

1) If $f_1(\xi)>0$ for $\xi\in (\xi_{i-1},\,\xi_i)$ and $f_1(\xi)>0$ for 
$\xi\in (\xi_{j},\,\xi_{j+1})$, then we define a function $f_2$ as
$$ \displaystyle
 f_2(\xi) = 
\left\{
\begin{array}{ll}
\frac{\gamma_1-f(\xi_{i-1})}{\xi_{i}-\xi_{i-1}} (\xi - \xi_{i}) + \gamma_1,  & \xi\in[\xi_{i-1},\,\xi_{i}], \\
\gamma_1,  & \xi\in[\xi_{i},\,\xi_{j}], \\
\frac{f(\xi_{j+1})-\gamma_1}{\xi_{j+1}-\xi_{j}} (\xi - \xi_{j}) + \gamma_1,  & \xi\in[\xi_{j},\,\xi_{j+1}], \\
\end{array}
\right.
$$
where
$\gamma_1:=\min\{\varepsilon/12,\,f(\xi_{i-1}),\,f(\xi_{j+1})\}.$ 

2) If $f_1(\xi)<0$ for $\xi\in (\xi_{i-1},\,\xi_i)$ and $f_1(\xi)<0$ for 
$\xi\in (\xi_{j},\,\xi_{j+1})$, then analogously to 1) we define 
$$ \displaystyle
 f_2(\xi) = 
\left\{
\begin{array}{ll}
\frac{\gamma_2-f(\xi_{i-1})}{\xi_{i}-\xi_{i-1}} (\xi - \xi_{i}) + \gamma_2,  & \xi\in[\xi_{i-1},\,\xi_{i}], \\
\gamma_2,  & \xi\in[\xi_{i},\,\xi_{j}], \\
\frac{f(\xi_{j+1})-\gamma_1}{\xi_{j+1}-\xi_{j}} (\xi - \xi_{j}) + \gamma_2,  & \xi\in[\xi_{j},\,\xi_{j+1}], \\
\end{array}
\right.
$$
where
$\gamma_2:= \max\{-\varepsilon/12,\,f(\xi_{i-1}),\,f(\xi_{j+1})\}.$

3) If the signs of $f_1$ are different on the intervals 
$(\xi_{i-1},\,\xi_i)$ and $(\xi_{j},\,\xi_{j+1})$, say, $f_1$ is positive on 
$(\xi_{i-1},\,\xi_i)$ and $f_1$ is negative on $(\xi_{j},\,\xi_{j+1})$, then we define $f_2$ as 
$$ \displaystyle
 f_2(\xi) = 
\left\{
\begin{array}{ll}
\frac{\gamma_3-f(\xi_{i-1})}{\xi_{i}-\xi_{i-1}} (\xi - \xi_{i}) + \gamma_3,  & \xi\in[\xi_{i-1},\,\xi_{i}], \\
\frac{\gamma_3-\gamma_4}{\xi_{j}-\xi_{i}} (\xi - \xi_{i}) + \gamma_3,  & \xi\in[\xi_{i},\,\xi_{j}], \\
\frac{f(\xi_{j+1})-\gamma_4}{\xi_{j+1}-\xi_{j}} (\xi - \xi_{j}) + \gamma_4,  & \xi\in[\xi_{j},\,\xi_{j+1}], \\
\end{array}
\right.
$$
where $\gamma_3:=\min\{\varepsilon/12,\,f(\xi_{i-1})\},$ $\gamma_4:=\max\{-\varepsilon/12,\,f(\xi_{j+1})\}.$
On the remaining part of $\mathbb{T}$ the functions $f_1$ and $f_2$ coincide.

 Thus, we obtain the piecewise linear function $f_2$ that has only finite number of roots and $\|f-f_2\|_C<\varepsilon /3.$ We need to keep the inequality 
$|f_2(\xi_k)|^2+|f_2(\xi_k+\pi)|^2\le 1,$ $k=0,\dots, 2n-1.$ The definition of $f_2$ implies that we modify the value of the function $f_1$ at points $\xi_i, \dots, \xi_j.$ To keep the inequality we change a bit (to keep the approximation $\|f-f_2\|_C<\varepsilon /3$ as well) the values of the function $f_1$ an the corresponding points $\xi_{N(i)}, \dots, \xi_{N(j)},$ where $N(i) = i+n,$ as $i=0,\dots, n-1,$ and $N(i)=i-n,$ as $i=n,\dots,2n-1.$ 

Finally, we remove pairs of symmetric roots and cycles of $f_2$, if any. Let $\tilde{\xi}$ be one of the symmetric roots or one of the roots from a cycle and $\tilde{\xi}\in[\xi_k,\,\xi_{k+1}).$ Then we replace the node of the polyline $(\xi_k,f_2(\xi_k))$ by the node $(\xi_k,f_2(\xi_k) \pm \varepsilon/12),$ where we choose the sign $``+''$ if $f(\xi_k)>0$ or $f(\xi_k)=0$ and $f(\xi_{k+1})>0$. For other cases we choose the sign $``-''$. 
 As a result, we obtain a new piecewise linear function
 $f_3$ such that $\|f-f_3\|_C< \varepsilon/2$. This function has only finite number of roots. It has  neither nontrivial cycles nor pairs of symmetric roots on $\mathbb{T}$. As on the previous step, to keep the inequality 
$|f_3(\xi_k)|^2+|f_3(\xi_k+\pi)|^2\le 1,$ $k=0,\dots, 2n-1,$ we modify the values of the function $f_2$ at the the corresponding points $\xi_{N(k)}$, if necessary.         

2. We need to check that 
\begin{equation}
\label{f3ineq}
|f_3(\xi)|^2+|f_3(\xi+\pi)|^2 \le 1
\end{equation}
 for any $\xi\in\mathbb{R}.$ Since $f_3$ and $f_3(\cdot+\pi)$ are linear on $[\xi_k,\xi_{k+1}],$ it follows that  
$$
\displaystyle \left(|f_3(\xi)|^2+|f_3(\xi+\pi)|^2\right)'' = 2(f'_3(\xi))^2+2(f'_3(\xi+\pi))^2 \ge 0,
$$ 
therefore the function $|f_3(\xi)|^2+|f_3(\xi+\pi)|^2$ is convex on $[\xi_k,\xi_{k+1}]$. The desired inequality for any $\xi \in [\xi_k,\xi_{k+1}]$ follows from the inequalities   $|f_3(\xi_k)|^2+|f_3(\xi_k+\pi)|^2 \le 1$ and $|f_3(\xi_{k+1})|^2+|f_3(\xi_{k+1}+\pi)|^2 \le 1.$ 

3. Let $\xi_k^j = \pi k/j,$ $k\in \mathbb{Z}.$ The interpolating trigonometric polynomial $H_j$ of the order $2j$ (the coefficient of $\cos 2j \xi$ is equal to zero) satisfying the conditions $H_j(\xi^j_k) = f_3(\xi^j_k),$ $H'_j(\xi^j_k) = 0$, $k=0,\dots, 2j-1,$ is defined by
$$
H_j(\xi) = \sum_{k=0}^{2j-1}f_3(\xi^j_k) t^{j}_k(\xi), 
\quad \mbox{ where } \quad
t^{j}_k(\xi) = \left(\frac{\sin j \xi}{2j \sin\frac{\xi -\xi^j_k}{2}}\right)^2.
$$
It is known (see \cite[page 51]{Tur}) that 
$\displaystyle 
\sum_{k=0}^{2j-1} t^{j}_k(\xi) =1,
$
and 
$\displaystyle 
\|H_j-f\|_{C} \rightarrow 0
$
as $j\to \infty.$ Let us check that $|H_j(\xi)|^2+|H_j(\xi+\pi)|^2 \le 1$ for any $\xi\in\mathbb{R}.$
Since $\displaystyle t^{j}_k(\xi+\pi) = \left(\frac{\sin j \xi}{2j \sin\frac{\xi -\left(\xi^j_k +\pi\right)}{2}}\right)^2 =t^{j}_{k+j}(\xi),$ it follows that 
$$
|H_j(\xi)|^2+|H_j(\xi+\pi)|^2 =
\left(\sum_{k=0}^{2j-1}f_3(\xi^j_k) t^{j}_k(\xi)\right)^2 +
\left(\sum_{k=0}^{2j-1}f_3(\xi^j_k) t^{j}_{k+j}(\xi)\right)^2 .
$$
Estimating $\displaystyle f_3(\xi^j_k)$ by $\displaystyle\left|f_3(\xi^j_k)\right|$ and changing the index of summation in the second sum, we obtain
$$
|H_j(\xi)|^2+|H_j(\xi+\pi)|^2 \le 
\left(\sum_{k=0}^{2j-1} \left|f_3(\xi^j_k)\right| t^{j}_k(\xi)\right)^2 +
\left(\sum_{k=0}^{2j-1} \left|f_3(\xi^j_{k+j})\right| t^{j}_{k}(\xi)\right)^2
$$
$$
=\sum_{k=0}^{2j-1} \left(\left|f_3(\xi^j_k)\right|^2 +\left|f_3(\xi^j_{k+j})\right|^2\right) \left(t^{j}_k(\xi)\right)^2
$$
$$
+2 \sum_{k<n} \left(\left|f_3(\xi^j_k)\right| \left|f_3(\xi^j_n)\right|+
\left|f_3(\xi^j_{k+j})\right| \left|f_3(\xi^j_{n+j})\right|\right)
t^{j}_k(\xi)
 t^{j}_n(\xi).
$$
Denote $a_k := \left|f_3(\xi^j_{k})\right|,$ 
by (\ref{f3ineq}), $a^2_k+a^2_{k+j}\le 1,$ so 
$$
a_k a_n +a_{k+j} a_{n+j} \leq a_k a_n +\sqrt{1-a_k^2} \sqrt{1-a_n^2} \le 
\max_{a_k,a_n \in [0,1]} \left(a_k a_n +\sqrt{1-a_k^2} \sqrt{1-a_n^2}\right) =1.
$$   
Therefore, taking into account that by (\ref{f3ineq}) 
$\displaystyle
\left|f_3(\xi^j_k)\right|^2 +\left|f_3(\xi^j_{k+j})\right|^2 \le 1,
$
we finally obtain 
$$
|H_j(\xi)|^2+|H_j(\xi+\pi)|^2 \le 
\sum_{k=0}^{2j-1}  \left(t^{j}_k(\xi)\right)^2
+2 \sum_{k<n} 
t^{j}_k(\xi)
 t^{j}_n(\xi)= \left(\sum_{k=0}^{2j-1}  t^{j}_k(\xi)\right)^2 =1.
$$

Now we claim that for a small enough $\varepsilon_1$ the polynomial $H_j$ has neither nontrivial cycles nor pairs of symmetric roots on $\mathbb{T}$. Indeed, let $f_3$ be defined explicitly as $f_3(\xi)=a_k \xi+b_k$ for $\xi \in [\xi_k,\xi_{k+1}],$ $k=0,\dots,2n-1,$ $\xi_0=0,$ $\xi_{2n}=2\pi.$ 
 We recall that  $\tilde{\xi}$ is a root and  we replace it to the new root $\tilde{\xi}'$ to remove pairs of symmetric roots and cycles of the function $f_2$. We denote by  $\alpha$ the minimum of $|\tilde{\xi}-\tilde{\xi}'|/2$ over all the pairs of the old and new roots $(\tilde{\xi},\,\tilde{\xi}')$  and  $a:=\min\{|a_k|:a_k \neq 0\}.$  Suppose the function $f_3$ has a root $\tilde{\xi}$ on the segment $[\xi_k,\xi_{k+1}]$. If $\|f_3-H_j\|_C<\varepsilon_1$, then 
$a_k \xi+b_k-\varepsilon_1 \leq H_j(\xi) \leq a_k \xi+b_k+\varepsilon_1$. In other words,  the plot of the polynomial $H_j$ lies inside the parallelogram bounded by the lines $y(\xi)=a_k \xi+b_k \pm \varepsilon_1$, $\xi=\xi_k$, $\xi=\xi_{k+1}.$ 
Therefore, the roots of $H_j$ can only be in the neighborhood of $\tilde{\xi},$ namely in the interval of the length $2\varepsilon_1/a_k$. If we provide the inequality 
$\varepsilon_1/a_k<\alpha$ for all the segments  $[\xi_k,\xi_{k+1}]$ that contains roots of $f_3$, then  it does mean that $H_j$ has neither nontrivial cycles nor pairs of symmetric roots on $\mathbb{T}$. To provide the inequality it is sufficient   to choose $\varepsilon_1< a \alpha.$

4. Since $H_j$ is a trigonometric polynomial, $H_j(0)=1,$ $|H_j(\xi)|^2+|H_j(\xi+\pi)|^2 \le 1$ for any $\xi\in\mathbb{R}$, by Remark \ref{rem1} it follows that $H_j$ is a refinement mask of the UEP. \hfill $\Box$

\begin{rem} It is easy to find the low-pass filter 
$\left(c^j_m(f)\right)_{m=-2j+1}^{2j-1}$ 
corresponding the refinement mask $H_j.$ 
$$
H_j(\xi) = \sum_{k=0}^{2j-1}f_3(\xi^j_k) t^{j}_k(\xi) = 
 \sum_{k=0}^{2j-1}f_3(\xi^j_k) \frac{1}{4 j^2} \left(\frac{\sin j (\xi - \xi^j_k)}{\sin((\xi-\xi^j_k)/2)}\right)^2 
$$
Applying the formula for the Fejer kernel, 
$$
\frac{1}{4 j} \left(\frac{\sin j x}{\sin(x/2)}\right)^2
 = \frac{1}{2} +\sum_{m=1}^{2j-1}\left(1-\frac{m}{2j}\right)\cos mx
= \sum_{m=-2j+1}^{2j-1}\frac{1}{2}\left(1-\frac{|m|}{2j}\right) e^{imx},
$$
we obtain 
$$
H_j(\xi) = \sum_{m=-2j+1}^{2j-1}\frac{1}{2j} \left(1-\frac{|m|}{2j}\right) \sum_{k=0}^{2j-1} f_3(\xi^j_k) e^{i m (\xi - \xi^j_k)}.
$$
Therefore, 
$$
c_m^j(f) = \frac{1}{2j} \left(1-\frac{|m|}{2j}\right) \sum_{k=0}^{2j-1} f_3(\xi^j_k) e^{-i m \xi^j_k}.
$$
Taking into account that $\xi^j_k = \pi k/j,$ we finally get 
$$
c_m^j(f) = \left(1-\frac{|m|}{2j}\right) \tilde{f}_3(m), 
$$
where $\displaystyle  \tilde{f}_3(m)=
 \frac{1}{2j}  \sum_{k=0}^{2j-1} f_3(\xi^j_k) e^{-i m \xi^j_k}$
is the inverse discrete Fourier transform of 
$\left(f_3(\xi^j_k)\right)_{k=0}^{2j-1}.$

\end{rem}




\end{document}